\documentclass[12pt, reqno]{amsart}

\usepackage[a4paper]{geometry}
\usepackage[matrix,arrow,curve,cmtip]{xy}

\usepackage{txfonts}

\usepackage{amssymb,amsmath,amsthm}

\def\identity{\operatorname{id}}

\def\End{\operatorname{End}}
\def\Ker{\operatorname{Ker}}
\def\Vir{\operatorname{Vir}}
\def\oo#1{\mathrel{{}_{(#1)}}}
\def\Cur{\mathrm{Cur}\,}
\def\SpanL{\mathrm{Span}\,}
\def\Res{\mathrm{Res}}
\newtheorem{theorem}{Theorem}
\newtheorem{proposition}{Proposition}
\newtheorem{lemma}{Lemma}
\newtheorem{corollary}{Corollary}
\theoremstyle{definition}
\newtheorem{definition}{Definition}
\newtheorem{example}{Example}

\begin{document}

\title{Homogeneous averaging operators on simple finite conformal Lie algebras} 

\author{Pavel Kolesnikov}
\email{pavelsk@math.nsc.ru}
\thanks{Supported by Russian Science Foundation (project 14-21-00065)}
\address{Novosibirsk State University}
\address{Sobolev Institute of Mathematics, Akad. Koptyug ave., 4, 630090 Novosibirsk, Russia}

\begin{abstract}
We describe all homogeneous averaging operators on finite simple conformal Lie algebras. 
In the case of current algebras, these operators are closely related to 
rational solutions of the classical Yang---Baxter equation.
\end{abstract}

\keywords{Conformal algebra, averaging operator, Yang---Baxter equation}

\maketitle 

\bibliographystyle{plain}

\section{Introduction}

In this note we show a close relation between 
classical Yang---Baxter equation (CYBE),
conformal algebras, 
and averaging operators on Lie algebras. 

Averaging operators initially appear in Reynolds' works in turbulence theory and later found many applications 
in analysis. We refer a reader to the well-written introduction 
of~\cite{GuoPei2014}
 for more details. 
For a Lie algebra $\mathfrak g$, an averaging operator is a linear map $T:\mathfrak g\to \mathfrak g$ 
such that $T([T(x),y])=[T(x),T(y)]$ for all $x,y\in \mathfrak g$. 

Conformal algebras were introduced by V.~Kac \cite{Kac1996}
 as ``Lie analogues'' of vertex algebras. 
Namely, the commutator of two chiral fields in 2-dimensional conformal field theory 
may be expressed in terms of coefficients taken from 
the singular part of their operator product expansion (OPE). Algebraic formalization of this expression 
leads to the notion of a (Lie) conformal algebra, a linear space equipped with a linear map $\partial $
and with a family of bilinear products $[\cdot \oo{n} \cdot ]$, where $n$ ranges over 
the set $\mathbb Z_+$ of nonnegative integers. 
It is natural to define an averaging operator on a conformal algebra $C$ as a $\partial $-invariant map $T$ satisfying 
$T([T(x)\oo{n} y]) = [T(x)\oo{n} T(y)]$ for all $x,y \in C$, $n\in \mathbb Z_+$.

One of the simplest examples of a conformal algebra 
is the current algebra $\Cur\mathfrak g=\mathbb C[\partial ]\otimes \mathfrak g$ 
over a Lie algebra~$\mathfrak g$. 
Other examples include Virasoro conformal algebra $\Vir $ and semi-direct product $\Vir \ltimes \Cur\mathfrak g$
\cite{Kac1996}.

An averaging operator on $\Cur \mathfrak g$ may be expressed 
in terms of a series of linear operators $T_n:\mathfrak g\to \mathfrak g$ in the following way:
$T(1\otimes a) = \sum_{n\ge 0} (-\partial )^n/n! \otimes T_n(a)$, $a\in \mathfrak g$. 
This leads to a map $T_\lambda : \mathfrak g\to \mathfrak g[\lambda ]$
given by $T_\lambda (a) = T(1\otimes a)|_{\partial =-\lambda }$, where $\lambda $
is a formal variable. 
The map $T_\lambda $ obtained satisfies 
\[
 T_{\lambda+\mu}([T_\lambda (a),b]) = [T_\lambda (a), T_\mu (b)], \quad a,b\in \mathfrak g.
\]
It turns out that the same equation holds for the singular part of every solution of the 
classical Yang---Baxter equation (CYBE) on the Lie algebra $\mathfrak g$.
Therefore, every solution of CYBE on a finite-dimensional Lie algebra $\mathfrak g$
determines an averaging operator on the corresponding current conformal algebra $\Cur \mathfrak g$.

The problem of description of all averaging operators on a conformal algebra (as well as on an ordinary algebra) 
seems to be hard. However, we can describe all {\em homogeneous} operators. An averaging
operator $T$ on current conformal algebra $\Cur \mathfrak g$ is said to be homogeneous 
if the collection of coefficients of polynomials from $T_\lambda (\mathfrak g)$ includes 
a Cartan subalgebra of $\mathfrak g$.
These results were announced in \cite{Kol2015AL}. 
We also describe all averaging operators 
on the Virasoro conformal algebra and, in the final section, all (not only homogeneous)
averaging operators on $\Cur\mathfrak{sl}_2$ up to conjugation with an automorphism of $\mathfrak{sl}_2$. 
The latter approach may also be applied to classification of Rota---Baxter operators 
\cite{BaiNiGuo2014}: 
we present a different proof of the main technical result of 
\cite{BaiNiGuo2014}
from this point of view.

\section{Averaging operators}

\begin{definition}
 Suppose $A$ is a linear algebra (not necessarily associative) over $\mathbb C$ with a product $\cdot : A\otimes A\to A$, and let $T$ be a linear 
transformation of the space $A$. Then $T$ is an {\em averaging operator} if 
\begin{equation}\label{eq:Averaging}
 T(T(a)\cdot b) = T(a)\cdot T(b) = T(a\cdot T(b))
\end{equation}
for all $a,b\in A$.
\end{definition}
If $A$ is commutative or anti-commutative then it is enough to check only one of the equations in~\eqref{eq:Averaging}.

Let us state here some simple examples.

\begin{example}\label{exmp:NonDegenerateCentroid}
Let $A$ be an algebra.
\begin{itemize}
 \item
If $\Gamma $ is a finite group of automorphisms of $A$ then 
the following transformation is an averaging operator on $A$:
\[
 T=T_\Gamma : a \mapsto \sum\limits_{g\in \Gamma} g(a), \quad a\in A.
\]

\item
If $T$ is an averaging operator on $A$ and $\sigma $ is an automorphism of $A$ 
then $T^\sigma = \sigma^{-1} T\sigma $ is also an averaging operator on $A$.

\item
 Let $T_1$ and $T_2$ be commuting averaging operators on $A$. 
Then $T_1T_2$ is an averaging operator on $A$. In particular, a power
of an averaging operator is also an averaging operator.

\item
 A surjective averaging operator $T$ on $A$ belongs to the centroid of $A$, i.e., 
such $T$ commutes with all operators of left and right multiplication in $A$.
In particular, if $A$ is a finite-dimensional simple algebra 
then a non-degenerate averaging operator $T$ is just a scalar map, 
$T=\alpha \identity$, $\alpha \in \mathbb C$.

\item
If $A$ is a finite-dimensional classically semisimple algebra, i.e., $A=\bigoplus\limits_{i=1}^n A_i$, where $A_i$ are simple algebras, 
then a non-degenerate averaging operator $T$ acts on $A$ as a direct sum of scalar maps on $A_i$, $i=1,\dots , n$. 
\end{itemize}
\end{example}

Averaging operators also occur in the theory of Leibniz algebras, the most popular and well-studied 
noncommutative analogues of Lie algebras. Namely, if $\mathfrak g$ is a Lie algebra with a product 
$[\cdot ,\cdot ]$ equipped with an averaging operator $T$ then the same space relative to new binary operation
\[
 \{a,b\}_T = [T(a),b],\quad a,b\in \mathfrak g,
\]
satisfies the Jacobi identity $\{x,\{y,z\}_T\}_T - \{y,\{x,z\}_T\}_T = \{\{x,y\}_T,z\}_T$, 
i.e., $\mathfrak g_T = (\mathfrak g, \{\cdot ,\cdot \}_T)$
is a Leibniz algebra. Note that $\Ker T$ is an ideal 
of $\mathfrak g_T$, and $\mathfrak g_T/\Ker T$ is a Lie algebra.
Moreover, every Leibniz algebra may be embedded into an algebra of the form $\mathfrak g_T$
for an appropriate Lie algebra $\mathfrak g$ with an averaging operator $T$ 
\cite{GubKol2014}.  

There is a sort of duality between averaging operators and Rota---Baxter operators. 
Recall that a Rota---Baxter operator of weight zero on a (non-associative, in general) algebra $A$ 
is a linear map $R:A\to A$ such that
\begin{equation}\label{eq:RotaBaxer}
 [R(x),R(y)] = R([x,R(y)]) + R([x,R(y)]),  \quad x,y\in A.
\end{equation}
Suppose $\mathfrak V$ is a variety of algebras (e.g., associative, alternative, Jordan), and $A\in \mathfrak V$.
If $T$ is an averaging operator on $A$ then the same linear space $A$ equipped with 
two new operations $x\vdash y = T(x)y$ and $x\dashv y = xT(y)$ is a system 
from the variety di-$\mathfrak V$ of $\mathfrak V$-dialgebras 
\cite{GubKol2011}. 
If $R$ is a Rota---Baxter operator of weight zero on $A$ then two new 
operations $x\succ y = R(x)y$, $x\prec y = xR(y)$ turn $A$ into a system from the variety 
pre-$\mathfrak V$ of dendriform $\mathfrak V$-dialgebras 
\cite{BaiGuoNi2012}.
A general statement for di-algebras 
\cite{GubKol2011} says that if $\mathfrak V$
is governed by a quadratic binary operad then 
\[
 (\text{di-}\mathfrak V)^! = \text{pre-}(\mathfrak V^!)
\]
where $(\cdot)^!$ denotes Koszul duality of operads.

\section{Conformal algebras}

\begin{definition}[{\cite{Kac1996}}]
A linear space $C$ equipped with a
linear operator $\partial : C\to C$ and with a family of bilinear ``products'' 
$[\cdot \oo{n} \cdot ]$,  $n\in \mathbb Z_+$, is called a {\em conformal algebra}
if for every $a,b\in C$ the following statements hold:
\begin{itemize}
 \item [(C1)] There exists $N=N(a,b)$ such that  $[a\oo{n} b] = 0$ for  $n\ge N$;
 \item[(C2)] $[\partial a \oo{n} b] = -n [a\oo{n-1} b]$;
 \item[(C3)] $[a\oo{n} \partial b] = \partial ([a\oo{n} b]) + n[a\oo{n-1} b]$.
\end{itemize}
If, in addition, 
\begin{equation}\label{eq:ConfAComm}
 [a\oo{n} b] = \sum\limits_{s\ge 0} (-1)^{n+s+1} \frac{1}{s!}\partial^{s}([b\oo{n+s} a])
\end{equation}
and 
\begin{equation}\label{eq:ConfJacobi}
 [a\oo{n} [b\oo{m} c]] - [b\oo{m} [a\oo{n} c]] = \sum\limits_{s\ge 0} \binom{n}{s} [[a\oo{n-s} b] \oo{m+s} c]
\end{equation}
for all $a,b,c\in C$, $n,m\in \mathbb Z_+$ then 
$C$ is said to be a Lie conformal algebra (also known as 
{\em Lie vertex algebra} \cite{FBZhvi2001}). 
\end{definition}

The following example of a conformal algebra, though simplest possible, 
is essential for our needs.

\begin{example}\label{exmp:Current}
 Suppose $\mathfrak g$ is a Lie algebra. Then 
$C=\mathbb C[\partial ]\otimes \mathfrak g$ 
equipped with 
\[
[ (1\otimes a)\oo{n} (1\otimes b)] = \delta_{n,0}\otimes [a,b],\quad a,b\in \mathfrak g
\]
(one may use (C2) and (C3) to expand these operations to the entire $C$)
turns into a Lie conformal algebra denoted by $\Cur \mathfrak g$
(current conformal algebra).
\end{example}

There also exists a non-trivial conformal Lie algebra of rank one.

\begin{example}\label{exmp:Virasoro}
 Consider a 1-generated free
$\mathbb C[\partial ]$-module 
  $C=\mathbb C[\partial ]v$  equipped with 
\[
[ v\oo{n}v ] = \delta_{n,0}\partial v + 2\delta_{n,1}v.
\]
This is a Lie conformal algebra called {\em Virasoro conformal algebra} $\Vir$.
\end{example}

It was established in \cite{DK1998} that every {\em finite} (i.e., finitely generated 
as a $\mathbb C[\partial]$-module) simple Lie conformal algebra is isomorphic 
either to $\Vir $ or to $\Cur \mathfrak g$ for a simple finite-dimensional Lie algebra~$\mathfrak g$.

Every conformal algebra $C$ gives rise to {\em annihilation algebra} $\mathcal A(C)$ constructed as follows 
\cite{BDK2001}.
Consider 
the space of formal power series
$\mathbb C[[t]]$ as a right $\mathbb C[\partial]$-module assuming $f(t)\partial = -f'(t)$. Then the linear space 
$\mathcal A(C)=\mathbb C[[t]]\otimes _{\mathbb C[\partial]} C$ may be equipped with a bilinear product defined by
\[
 (t^n\otimes _{\mathbb C[\partial]} a)(t^m\otimes _{\mathbb C[\partial]} b)
 = 
\sum\limits_{s\ge 0}  \binom{n}{s} t^{n+m-s} \otimes_ {\mathbb C[\partial]} [a\oo{s} b] ,
\quad a,b\in C.
\]
This operation is well-defined and may be extended on the entire $\mathcal A(C)$ by continuity (with respect to $t$-adic 
topology).
If $C$ is a Lie conformal algebra then $\mathcal A(C)$ is a Lie algebra.  The converse is true for torsion-free conformal algebras, see
\cite[Section 11.3]{BDK2001}.

Suppose $C$ is a finite Lie conformal algebra,  $T:C\to C$ is a $\mathbb C[\partial ]$-linear map 
such that 
\begin{equation}\label{eq:AverConformal}
 T([T(a)\oo{n} b]) = [T(a)\oo{n} T(b)],\quad  a,b\in C,\ n\in \mathbb Z_+.
\end{equation}
Let us call such  a map by an {\em averaging operator} on the conformal algebra~$C$.

The map $C\mapsto \mathcal A(C)$ is a functor from the category of conformal Lie algebras to 
the category of topological differential Lie algebras. In particular, it is easy to see that if $T$ is an averaging operator on $C$ 
then $\mathcal A(T): t^n\otimes_{\mathbb C[\partial ]} a \mapsto t^n\otimes_{\mathbb C[\partial ]} T(a)$
is an averaging operator on $\mathcal A(C)$. Such an operator  necessarily commutes with the derivation $\partial_t=\partial/\partial t$
and continuous with respect to $t$-adic topology. Conversely, every $\partial_t$-invariant continuous averaging operators on $\mathcal A(C)$
gives rise to an averaging operator on the conformal algebra~$C$ via the reconstruction functor~\cite{BDK2001}.

The following statements may be easily shown by means of either direct routine computation, 
or categorical formalism 
\cite{BDK2001}, or using annihilation algebras. 
Suppose $C$ is a Lie conformal algebra and $T$ is an averaging operator on $C$.
Then the same module $C$ equipped with new family of operations 
$\{a\oo{n} b\}_T = [T(a)\oo{n} b]$ is a conformal algebra $C_T$ satisfying 
\eqref{eq:ConfJacobi}, i.e., it may be naturally called a {\em conformal Leibniz algebra}.
As in the case of ordinary Leibniz algebras, it is easy to see that 
$C_T/\Ker T $ is a Lie conformal algebra.

\section{Singular part of CYBE solution}

Let $\mathfrak g$ be a finite-dimensional Lie algebra. Suppose $X$ is a $(\mathfrak g\otimes \mathfrak g)$-valued 
function of a complex variable $u$ which is meromorphic at $u=0$, i.e., $X(u)$ is presented by a Laurent series 
in a neighborhood $\mathcal D\subset \mathbb C$ of the origin.
As usual, if $X=\sum\limits_i x_i'\otimes x_i''\in \mathfrak g\otimes \mathfrak g$ then $X^{12}$
stands for $\sum\limits_i x_i'\otimes x_i''\otimes 1\in U(\mathfrak g)^{\otimes 3}$, and the same 
convention determines $X^{13}$, $X^{23}$.

The classical Yang---Baxter equation (CYBE) is the functional equation 
\begin{equation}\label{eq:CYBE_tensor}
 [X^{12}(u), X^{13}(u+v)] + [X^{12}(u), X^{23}(v)]+[X^{13}(u+v),X^{23}(v)] =0.
\end{equation}
Solutions of CYBE are of great interest for pure algebra since they are related to quantizations 
of Lie bialgebras, see, e.g., \cite{MontStolinZelmanov}
 and references therein.
Constant solutions of CYBE and its generalizations are 
 related with Rota---Baxter operators and their generalizations 
\cite{BaiNiGuo2014}.
A.~Belavin and V.~Drinfeld 
\cite{BelDr} classified non-degenerate solutions of CYBE. 

If $\mathfrak g$ is a semisimple algebra then its Killing form $\langle \cdot , \cdot \rangle $
is non-degenerate, and one may identify $\mathfrak g\otimes \mathfrak g$ with $\End \mathfrak g$
by the natural rule
\[
 a\otimes b \mapsto \varphi_{a\otimes b}, \quad \varphi_{a\otimes b}(x) = \langle a,x\rangle b,
\]
for $a,b,x\in \mathfrak g$.
For example, the Casimir tensor $\Omega \in \mathfrak g\otimes \mathfrak g $ corresponds to the identity map. It is easy to see that 
$\varphi_{b\otimes a} = \varphi_{a\otimes b}^*$, where $*$ stands for the 
conjugation in $\End \mathfrak g$ relative to the Killing form.

Therefore, a meromorphic tensor-valued 
function $X$ corresponds to a meromorphic operator-valued function $P_u = \varphi_{X(u)}$, $u\in \mathcal D$. 
Straightforward computation shows that \eqref{eq:CYBE_tensor} is equivalent to the following operator equation:
\begin{equation}\label{eq:CYBE_operator}
 P_{u+v}([x,P^*_u(y)]) - P_v([P_u(x),y])+[P_{u+v}(x),P_v(y)]=0
\end{equation}
for $x,y\in \mathfrak g$.
In particular, if $P_u=R$ is a skew-symmetric constant then  $R$ satisfies \eqref{eq:RotaBaxer},
i.e., $R$ is a Rota---Baxter operator of weight zero.

\begin{example}\label{exmp:AveToCYBE}
 Suppose $T$ is an averaging operator on a semisimple Lie algebra $\mathfrak g$
such that 
 $T([T^*(x),y] - [T(x),y])=0$ for every $x,y\in \mathfrak g$. 
Then 
\[
 P_u (a)= \dfrac{1}{u} T(a), \quad a\in \mathfrak g,
\]
is a solution of the operator CYBE~\eqref{eq:CYBE_operator}.

For example, $T$ may be symmetric relative to the Killing form.
\end{example}

It was established in 
\cite{BelDr} 
for a simple Lie algebra $\mathfrak g$ that if $P_u$ is a
non-degenerate solution of CYBE (that is, $\det P_u\ne 0$ at some point $u\in \mathcal D$)
then the singular part of $P_u$ is of the form $\dfrac{\lambda \identity}{u}$, $\lambda \in \mathbb C$.

In general, a solution of \eqref{eq:CYBE_operator} may have a more complicated singular part, 
which may be presented by its generating function---a polynomial in a formal variable $\lambda $:
\begin{equation}\label{eq:P_Residue}
 T_\lambda = \Res_{u=0} P_u \exp(\lambda u) \in (\End \mathfrak g)[\lambda ].
\end{equation}
Indeed, if 
\[
 P_u = \frac{1}{u^{N+1}} T_N + \dots + \frac{1}{u}T_0 + \dots ,\quad T_k\in \End \mathfrak g, 
\]
then 
$T_\lambda  = \lambda^{(N)}T_N + \dots + T_0$,
where $\lambda^{(n)}$ stands for $\lambda^n/n!$.

Let us deduce an equation on $T_\lambda $.

\begin{theorem}\label{thm:CYBE_to_ConfAve}
 Let $\mathfrak g$ be a finite-dimensional Lie algebra and let $P_u$ be a meromorphic  solution 
of CYBE \eqref{eq:CYBE_operator}. Then $T_\lambda $ given by \eqref{eq:P_Residue}
satisfies 
\begin{equation}\label{eq:Conformal_Average}
 T_{\lambda+\mu}([T_\lambda (x),y]) = [T_\lambda (x), T_\mu(y)] \in \mathfrak g[\lambda, \mu ]
\end{equation}
for all $x,y\in \mathfrak g$.
\end{theorem}

\begin{proof}
Let us multiply \eqref{eq:CYBE_operator} by $\exp((\lambda +\mu)v)\exp (\lambda u)$
and integrate the expression obtained:
\begin{multline}\label{eq:CYBE_integral}
\oint\limits_{|u|=r}\oint\limits_{|v|=d}
\exp((\lambda +\mu)v)\exp (\lambda u) \left( P_v([P_u(x),y])- P_{u+v}([x,P^*_u(y)]) \right )dvdu \\
=
 \oint\limits_{|u|=r}\oint\limits_{|v|=d}
\exp(\lambda(u+v) \exp(\mu v)
[P_{u+v}(x),P_v(y)] \,dvdu
\end{multline}
where $d<r$ and $d$ is small enough for the entire circle $|u|\le 2r$ to lie in $\mathcal D$.
The first summand of the left-hand side of \eqref{eq:CYBE_integral} 
obviously provides $T_{\lambda+\mu }([T_\lambda (x),y])$, 
the second one is equal to zero 
since for every fixed $u$, $|u|=r$, 
$P_{u+v}([x,P^*_u(y)])$ considered as a function on $v$ has no poles in $|v|\le d<r$.
The right-hand side of \eqref{eq:CYBE_integral}  can be calculated by 
substitution:
\begin{multline}\nonumber
\oint\limits_{|u|=r}\oint\limits_{|v|=d}
 \exp((\lambda +\mu)v)\exp (\lambda u)[P_{u+v}(x),P_v(y)] \, dvdu \\
=
\oint\limits_{|v|=d}\oint\limits_{|u|=r}
 [\exp (\lambda(u+v))P_{u+v}(x),\exp(\mu v)P_v(y)] \,dudv \\
=
\oint\limits_{|v|=d}\oint\limits_{\substack{w=v+u \\ |u|=r }}
 [\exp (\lambda w) P_w(x),\exp(\mu v)P_v(y)] \,dwdv ,
\end{multline}
which gives $[T_\lambda (x), T_\mu (y)]$
since $w=0$ is the only singular point in the interior of the curve $w=v+u$, $|u|=r$.
Therefore, \eqref{eq:Conformal_Average} holds.
\end{proof}

\begin{definition}
Let us call an operator $T_\lambda : \mathfrak g\to \mathfrak g[\lambda ]$
satisfying \eqref{eq:Conformal_Average} by {\em conformal averaging operator} on~$\mathfrak g$.
\end{definition}

The explanation of the term ``conformal'' comes from the following observation. 
Suppose $C=\Cur \mathfrak g$ is the current conformal algebra over $\mathfrak g$,
and let $T_\lambda $ be a conformal averaging operator on~$\mathfrak g$.
It turns out that $T_\lambda $ is a form of an averaging operator on 
the conformal algebra $\Cur \mathfrak g$.
 
\begin{proposition}\label{prop:ConformalAve-Cur}
Assume $T_\lambda $ is a conformal averaging operator on $\mathfrak g$,
\[
 T_\lambda (a) = \sum\limits_{n\ge 0} \lambda^{(n)} T_n(a), \quad a\in \mathfrak g,
\]
where $T_n$ are linear transformations of~$\mathfrak g$.
Consider $\mathbb C [\partial ]$-linear map
$T : \Cur \mathfrak g \to \Cur \mathfrak g$, 
given by 
\[
  T(1\otimes a) = T_{-\partial} (a) = \sum\limits_{n\ge 0} (-\partial) ^{(n)}\otimes T_n(a).
\]
Then $T$ satisfies \eqref{eq:AverConformal}.
\end{proposition}

\begin{proof}
It is enough to prove the statement for $x,y\in \mathfrak g\simeq 1\otimes \mathfrak g$.
Compare coefficients at $\lambda^{(n)}\mu^{(m)}$ in \eqref{eq:Conformal_Average}:
\begin{equation}\label{eq:ConfAssoc-1}
[T_n(x),T_m(y)] = \sum\limits_{t\ge 0} \binom{n}{t} T_{m+t}([T_{n-t}(x),y]).
\end{equation}
Conversely (replace $\mu$ with $\mu-\lambda$), 
\begin{equation}\label{eq:ConfAssoc-2}
T_m([T_n(x),y]) = \sum\limits_{s\ge 0} (-1)^s\binom{n}{s} [T_{n-s}(x),T_{m+s}(y)].
\end{equation}
It remains to calculate left- and right-hand sides of the desired relation in $\Cur\mathfrak g$:
\[
\begin{aligned}
 T([T(x)\oo{n} y]) & = \sum\limits_{m\ge 0} (-\partial )^{(m)}\otimes T_m([T_n(x),y]), \\
 [T(x)\oo{n}T(y)] & = \sum\limits_{k,s\ge 0}\binom{n}{k} [T_k(x)\oo{n-k} (-\partial )^{(s)}T_s(y)]  \\
 & = \sum\limits_{k+s=n}(-1)^s \binom{n}{k} [T_k(x),T_s(y)].
\end{aligned}
\]
\end{proof}

Therefore, every solution of CYBE \eqref{eq:CYBE_operator} gives rise to a conformal 
averaging operator $T_\lambda $ on $\mathfrak g$ and thus induces a Leibniz conformal algebra structure 
on $\mathbb C [\partial ]\otimes \mathfrak g$ which is $(\Cur\mathfrak g)_T$. 
Then $(\Cur\mathfrak g)_T/\Ker T$ is a Lie conformal algebra.
The purpose of the next section is to describe such Lie conformal algebras.

\section{Description of homogeneous conformal averaging operators}

Let $T_\lambda $ be a conformal averaging operator on a finite-dimensional Lie algebra $\mathfrak g$.
For a subspace $B$ of $\mathfrak g$ denote by 
$ T_*(B)$ the linear span of $T_\alpha (b)$, $b\in B$, $\alpha \in \mathbb C$.
Note that if $B$ satisfies the condition $[T_*(B),B]\subseteq B$
then $T_*(B)$ is a subalgebra of $\mathfrak g$ by \eqref{eq:ConfAssoc-1}.
In particular, $T_*(\mathfrak g)$, $T_*(T_*(\mathfrak g))$
are subalgebras of $\mathfrak g$.
Let us call $T_\lambda $ {\em non-degenerate} if $T_*(\mathfrak g)=\mathfrak g$.
If $\mathfrak g$ is semisimple and $T_*(\mathfrak g) $ contains a Cartan subalgebra $\mathfrak h$
of $\mathfrak g$ then $T_\lambda $ is called {\em homogeneous}.

\begin{lemma}\label{lem:Semisimple-NonDeg}
 Let $\mathfrak g$ be a finite-dimensional semisimple Lie algebra. 
Then every non-degenerate conformal averaging operator on $\mathfrak g$ 
is just an ordinary averaging operator.
\end{lemma}

\begin{proof}
 Suppose $T_\lambda = T_0 + \lambda T_1 +\dots + \lambda^{(N)}T_N$. 
It follows from \eqref{eq:ConfAssoc-1} that $T_N(\mathfrak g)$ is an ideal of $T_*(\mathfrak g)$, 
moreover, if $N>0$ then this ideal is abelian. Hence, $N=0$, and $T_\lambda(a)=T_0(a)$ for 
every $a\in \mathfrak g$, where $T_0$ is a non-degenerate averaging operator on $\mathfrak g$. 
All such operators are given by Example~\ref{exmp:NonDegenerateCentroid}.
\end{proof}

\begin{lemma}\label{lem:Reductive-NonDeg}
 Let $\mathfrak g$ be a finite-dimensional reductive Lie algebra, 
$\mathfrak g = \mathfrak g_0 \oplus Z$, where $\mathfrak g_0$ is semisimple, $Z$ is the center of $\mathfrak g$.
Suppose $T$ is a conformal averaging operator on $\mathfrak g$ such that 
$\mathfrak g_0\subseteq T_*(\mathfrak g)$.
Then $T_*(\mathfrak g_0)=\mathfrak g_0$, $T_*(Z)\subseteq Z$.
\end{lemma}

As a corollary, $T_\lambda $ is a non-degenerate conformal averaging operator on $\mathfrak g_0$
described by Lemma~\ref{lem:Semisimple-NonDeg}. Thus, an ordinary averaging 
operator on $\mathfrak g_0$ and an arbitrary map $Z\to Z[\lambda ]$,
define a conformal averaging operator $T_\lambda $ on~$\mathfrak g$.

\begin{proof}
Relation \eqref{eq:ConfAssoc-1} implies that $T_m(z)$ commutes with 
$\mathfrak g_0\subseteq T_*(\mathfrak g)$ for $z\in Z$.
Hence, $T_*(Z)\subseteq Z$, and  the induced map 
$\bar T_\lambda : \mathfrak g/Z\to \mathfrak g/Z[\lambda ]$ 
is a non-degenerate conformal averaging operator on $\mathfrak g_0\simeq  \mathfrak g/Z$.
By Lemma \ref{lem:Semisimple-NonDeg} $\bar T_\lambda = T_0: \mathfrak g_0\to \mathfrak g_0$.
Therefore, 
\[
 T_\lambda (a) = T_0(a) + \zeta_\lambda (a), \quad a\in \mathfrak g_0,
\]
where $\zeta_\lambda : \mathfrak g_0\to Z[\lambda ]$.

It follows from 
\eqref{eq:Conformal_Average} that
$ T_{\mu}([T_0(a),b]) = [T_0(a),T_\mu(b)]$, 
so $\zeta_{\mu}([T_0(a),b]) = 0$
for all $a,b\in \mathfrak g_0$. Hence, $\zeta_\lambda(\mathfrak g_0) =0$.
\end{proof}

If $\mathfrak g$ is finite-dimensional then $T_\lambda $ is 
defined by a finite family of linear operators 
$T_n:\mathfrak g\to \mathfrak g$, $n=0,\dots, N$,
as in Proposition~\ref{prop:ConformalAve-Cur}.

Consider the following chain of subspaces in $\mathfrak g$:
\[
 T_*(\mathfrak g) = T_{(0)}(\mathfrak g) \supseteq T_{(1)}(\mathfrak g)\supseteq \dots \supseteq T_{(N)}(\mathfrak g),
\]
where 
\[
 T_{(n)}(\mathfrak g) = T_{n}(\mathfrak g) + \dots + T_{N}(\mathfrak g), 
\quad 
n=0,1,\dots , N.
\]
Suppose $T_{(N+1)}(\mathfrak g)=\{0\}$.
Relations \eqref{eq:ConfAssoc-1} imply that 
$T_{(n)}(\mathfrak g)$ is an ideal of $T_*(\mathfrak g)$. 
Moreover, for every $n=0,1,\dots, N$ the induced map
$\bar T_n: \mathfrak g\to T_*(\mathfrak g)/T_{(n+1)}(\mathfrak g)$
given by  
$\bar T_n(a) = T_n(a)+T_{(n+1)}(\mathfrak g)$
is a homomorphism of $T_*(\mathfrak g)$-modules:
\begin{equation}\label{eq:QuotientAveraging}
 [T_k(x),T_n(a)] - T_n([T_k(x),a]) \in T_{(n+1)}(\mathfrak g),\quad x,a\in \mathfrak g.
\end{equation}

\begin{theorem}\label{thm:ConfAveraging_homog}
 Let $T_\lambda $ be a homogeneous conformal averaging operator on 
a semisimple finite-dimensional Lie algebra $\mathfrak g$.
Then $T_*(\mathfrak g)=\mathfrak g_0\oplus Z$ is a reductive Lie algebra
and $\mathfrak g_0\subseteq T_*(T_*(\mathfrak g))$.
\end{theorem}

\begin{proof}
Recall that $T_\lambda $ is homogeneous if $T_*(\mathfrak g)$
 contains a Cartan subalgebra~$\mathfrak h$ of~$\mathfrak g$.
Suppose $\Delta \subseteq \mathfrak h^*$ is the root system of $\mathfrak g$, and
\[
 \mathfrak g = \mathfrak h \oplus \sum\limits_{\alpha \in \Delta }\mathbb Cx_\alpha 
\]
is the root decomposition of $\mathfrak g$ relative to~$\mathfrak h$. 
Then $T_*(\mathfrak g)$ is a homogeneous subalgebra relative to this root grading.
Denote 
\[
 \Delta' = \{ \alpha \in \Delta \mid x_\alpha \in T_*(\mathfrak g) \}.
\]

To prove reductivity of $T_*(\mathfrak g)$ it is enough to show that 
$\Delta'$ is symmetric. If $\alpha \in \Delta'$ implies $-\alpha \in \Delta'$
then $\Delta'$ satisfies all necessary axioms of a root system, 
and the subalgebra $\mathfrak g_0\subseteq T_*(\mathfrak g)$
generated by $\{x_\alpha \mid \alpha \in \Delta'\}$ is a semisimple 
Lie algebra with Cartan subalgebra 
$\mathfrak h_0 = \SpanL\{h_\alpha =[x_\alpha ,x_{-\alpha }] \mid \alpha \in \Delta'\}$.
Finally, $\mathfrak h = \mathfrak h_0\oplus \mathfrak h_0^\perp$ relative to the Killing form, 
and $[\mathfrak h_0^\perp, x_\alpha ]=0$ for all $\alpha \in \Delta'$.
Hence, 
$T_*(\mathfrak g) = \mathfrak g_0 \oplus \mathfrak h_0^\perp$ 
is a reductive Lie algebra.

Let us show symmetry of $\Delta'$.
Assume $x_\alpha \in T_*(\mathfrak g)$. Choose the maximal $n\in \{0,\dots, N\}$ such that  
$x_\alpha \in T_{(n)}(\mathfrak g)$:
\[
 x_\alpha = T_n(y) + b, \quad y\in \mathfrak g, \ b\in T_{(n+1)}(\mathfrak g).
\]
Note that \eqref{eq:QuotientAveraging} implies 
\[
 T_n(x_\gamma )   \in \xi_{n,\gamma } x_\gamma +  T_{(n+1)}(\mathfrak g),
\]
for $\xi_{n,\gamma}\in \mathbb C$, $\gamma \in \Delta $,
and 
$T_n(h)\in \mathfrak h + T_{(n+1)}(\mathfrak g)$ for $h\in \mathfrak h$
since $T_\lambda $ is 
homogeneous.
Suppose 
$y = h + \zeta_\alpha x_\alpha +\sum\limits_{\beta \in \Delta\setminus \{\alpha \} } \zeta_\beta x_\beta $ is the root decomposition of $y$.
Then 
$T_{(n+1)}(\mathfrak g)\ni 
 x_\alpha -T_n(y) = 
(1-\xi_{n,\alpha }\zeta_{\alpha })x_\alpha -T_n(h)- \sum\limits_{\beta \in \Delta\setminus \{\alpha \} } \xi_{n,\beta}\zeta_\beta x_\beta$.
Since $T_{(n+1)}(\mathfrak g)$ is a homogeneous subalgebra relative to the root grading and 
$x_\alpha \notin T_{(n+1)}(\mathfrak g)$, we have $1-\xi_{n,\alpha }\zeta_\alpha =0$. 
Thus, $\xi_{n,\alpha }\ne 0$.

Assume $\xi_{n,-\alpha }$=0:
\[
 T_n(x_{-\alpha }) \in T_{(n+1)}(\mathfrak g).
\]
Then $[x_\alpha , T_n(x_{-\alpha })]\in T_{(n+1)}(\mathfrak g)$.
On the other hand, \eqref{eq:QuotientAveraging} implies 
\[
[x_\alpha , T_n(x_{-\alpha })] -  T_n([x_\alpha ,x_{-\alpha }])\in T_{(n+1)}(\mathfrak g),
\]
and thus $T_n(h_\alpha) = T_n([x_\alpha , x_{-\alpha }])\in T_{(n+1)}(\mathfrak g)$.
However, 
\[
 [x_\alpha , T_n(h_\alpha )] - T_n([x_\alpha , h_\alpha ])\in T_{(n+1)}(\mathfrak g)
\]
provides contradiction to $x_\alpha \notin T_{(n+1)}(\mathfrak g)$ since
$[h_\alpha , x_\alpha ]= (\alpha ,\alpha )x_\alpha \ne 0$.

Therefore, $T_n(x_{-\alpha }) \in \xi_{n,-\alpha} x_{-\alpha }+ T_{(n+1)}(\mathfrak g)$,
$\xi_{n,-\alpha }\ne 0$, and thus $x_{-\alpha }\in T_*(\mathfrak g)$. 

We have shown that if $x_\alpha \in T_*(\mathfrak g)$
then there exists $n$ such that 
$T_n(x_\alpha ) = \xi_{n,\alpha } x_\alpha + T_{(n+1)}(\mathfrak g)$, 
$x_\alpha \notin T_{(n+1)}(\mathfrak g)$. 
Therefore, every such $x_\alpha $ belongs to $ T_*(T_*(\mathfrak g))$. 
This implies $\mathfrak g_0\subseteq T_*(T_*(\mathfrak g))$.
\end{proof}

Theorem \ref{thm:ConfAveraging_homog} allows us to find explicit description 
of homogeneous conformal averaging operators on
a finite-dimensional semisimple Lie algebra $\mathfrak g$ (i.e., on $\Cur \mathfrak g$).

\begin{corollary}
 Under the conditions of Theorem \ref{thm:ConfAveraging_homog}, 
$\mathfrak g = T_*(\mathfrak g) \oplus \Ker T_\lambda $.
\end{corollary}

\begin{proof}
By  Theorem \ref{thm:ConfAveraging_homog},
 $T_\lambda $ is a conformal averaging operator on the reductive Lie algebra $T_*(\mathfrak g)$:
there  exists root subsystem $\Delta'\subseteq \Delta $ such that 
$T_*(\mathfrak g)$ is generated by $\mathfrak h$ and $x_\alpha $, $\alpha \in \Delta'$.
Moreover, $T_*(T_*(\mathfrak g))$ contains all $x_\alpha $ for $\alpha \in \Delta'$, so the semisimple 
factor $\mathfrak g_0$ of $T_*(\mathfrak g)$ lies in $T_*(T_*(\mathfrak g))$.
Lemma \ref{lem:Reductive-NonDeg} describes how $T_\lambda $ acts on $T_*(\mathfrak g)$, 
in particular, $T_*(\mathfrak h)\subseteq \mathfrak h$.

It remains to show that $T_\lambda (x_\beta )=0$ for all $\beta \in \Delta\setminus \Delta'$.
Assume there exists such a root $\beta $ that $T_\lambda (x_\beta )\ne 0$, $x_\beta \notin T_*(\mathfrak g)$. 
Let us choose maximal $n\ge 0$ such that $T_n(x_\beta )\ne 0$. 
Then \eqref{eq:QuotientAveraging} and $T_k(h)\in \mathfrak h$ (for all $k\ge 0$, $h\in \mathfrak h$)
imply 
$[h,T_n(x_\beta )] = \beta(h) T_n(x_\beta )$, i.e., $T_n(x_\beta )\in \mathbb C x_\beta $. 
As $x_\beta \notin T_*(\mathfrak g)$, we have $T_n(x_\beta )=0$, a contradiction.
\end{proof}

Finally, we may describe all homogeneous conformal averaging operators on 
a finite-dimensional semisimple Lie algebra $\mathfrak g$ 
with Cartan subalgebra $\mathfrak h$ and root system $\Delta $ as follows.

For a root subsystem $\Delta'$ of $\Delta $, 
\[
 T_\lambda (x_\alpha) = \xi_\alpha x_\alpha , \quad T_\lambda (h_\alpha )=\xi_\alpha h_\alpha ,
\]
where $\alpha \in \Delta'$, $h_\alpha = [x_\alpha , x_{-\alpha }]$. 
Note that $\xi_\alpha \in \mathbb C$ are nonzero constants that depend on the decomposition 
of the semisimple Lie algebra $\mathfrak g_0$ generated by $x_\alpha $, $\alpha \in \Delta'$, 
into simple summands. 
The subalgebra $\mathfrak h_0^\perp = \{h\in \mathfrak h \mid \alpha (h)=0, \ \alpha \in \Delta'\}$
is invariant with respect to $T_\lambda $, and there are no restrictions on 
linear functions 
$T_n|_{\mathfrak h_0^\perp}$ ($0\le n\le N$).
Other root spaces $\mathbb Cx_\beta $, $\beta \in \Delta\setminus \Delta'$, belong to the kernel of $T_\lambda $.

Straightforward computation shows 
\begin{equation}\label{eq:Complement-h}
[\mathfrak h_0^\perp , T_*(\mathfrak g)] =0, 
\quad 
[\mathfrak g, \mathfrak h_0^\perp ] \subseteq \Ker T_\lambda  
\end{equation}
since $[x_\alpha , h]=0$ for all $\alpha \in \Delta '$, $h\in \mathfrak h_0^\perp $.

Therefore, we have 

\begin{corollary}
For a homogeneous conformal averaging operator $T_\lambda $ on $\mathfrak g$
the induced Lie conformal algebra structure 
$(\Cur \mathfrak g)_T/\Ker T$ is given by split null extension $(\Cur \mathfrak g_0)\oplus (\Cur \mathfrak h_0^\perp)$.
\end{corollary}

It is worth mentioning that, in general, a conformal averaging operator $T_\lambda $ may not be a singular 
part of a solution of CYBE. However, it turns out that all homogeneous conformal averaging operators 
actually come from some meromorphic solutions of CYBE.

\begin{corollary}\label{cor:SingularCYBE_solution}
 Let $T_\lambda $ be a homogeneous conformal averaging operator on 
a finite-dimensional semisimple Lie algebra $\mathfrak g$.  
Then the operator-valued meromorphic function $P: \mathbb C\setminus\{0\} \to \End \mathfrak g $ given by
\[
 P_u(a) = \sum\limits_{n\ge 0} \dfrac{1}{u^{n+1} }T_n(a), \quad a\in \mathfrak g, \ u\in \mathbb C^*,
\]
is a solution of the classical Yang---Baxter equation \eqref{eq:CYBE_operator}.
\end{corollary}

\begin{proof}
It is easy to verify that 
\[
 T^*_\lambda  = \sum\limits_{n=0}^N \lambda ^{(n)} T^*_n,
\]
where $T_n^*$ is the conjugate of $T_n$ relative to the Killing form, 
is also a homogeneous conformal averaging operator with the same root subsystem $\Delta'$ as 
$T_\lambda $.
Indeed, the algebra $\mathfrak g$ splits into a direct sum of subspaces
$\mathfrak{g}= \mathfrak{g}_0\oplus K \oplus Z$, where $K=\Ker T_\lambda$, $Z=\mathfrak h_0^\perp$.
These subspaces are pairwise orthogonal relative to the Killing form, hence, $T_\lambda ^*$ is also 
a homogeneous conformal averaging operator, and 
$\Ker T_\lambda|_{\mathfrak{g}_0\oplus K} = \Ker T^*_\lambda|_{\mathfrak{g}_0\oplus K}$.
Therefore, it is enough to check whether~\eqref{eq:CYBE_operator} holds if either of $x$ or $y$ belong to $Z$.

If $x\in Z$ then 
 $P_u(x)\in Z$ for every $u\in \mathbb C\setminus \{0\}$
and thus all three summands of~\eqref{eq:CYBE_operator} are equal to zero due to~\eqref{eq:Complement-h}.
Finally, for $y\in Z$ equation \eqref{eq:CYBE_operator}
holds by the same reason.
\end{proof}

\section{Averaging operators on some conformal algebras}

In this section, we describe all (not necessarily homogeneous) 
averaging operators on the Virasoro conformal algebra $\Vir $ and on 
$\Cur \mathfrak {sl}_2$. The last case includes, in particular, all averaging operators on 
$\mathfrak{sl}_2$. 

\begin{proposition}\label{prop:VirAve}
 Let $T$ be an averaging operator on the Virasoro conformal algebra 
$\Vir = \mathbb C[\partial ]v$. 
Then $T(v)=\alpha v$, $\alpha \in \mathbb C$.
\end{proposition}
 
\begin{proof}
The useful notion of a $\lambda $-product 
\cite{DK1998} allows us to consider 
one generating series in a formal variable $\lambda $ instead of $n$-products in conformal algebras.
Namely, 
$[v_\lambda v] = [v\oo{0}v] +\lambda [v\oo{1}v] = (\partial +2\lambda )v$.
Suppose $T(v)=f(\partial )v$, then 
\[
 T([T(v)_\lambda v]) = [T(v)_\lambda T(v)]
\]
implies 
\[
 T(f(-\lambda )(\partial +2\lambda )v) = f(-\lambda )f(\partial +\lambda )(\partial+2\lambda )v.
\]
 If $f(\partial )\ne 0$ then $f(\partial )=f(\partial +\lambda )$, i.e., $f$ is a constant polynomial.
\end{proof}

Let us classify all conformal averaging operators $T_\lambda $ on $\mathfrak {sl}_2$, 
they correspond to averaging operators on the conformal algebra $\Cur\mathfrak{sl}_2$.
It is natural to perform the classification up to conjugation with an automorphism of 
$\mathfrak{sl}_2$ (see Example~\ref{exmp:NonDegenerateCentroid}).

Note that if $T_\lambda $ is a conformal averaging operator on $\mathfrak g$
and $\sigma $ is an automorphism of $\mathfrak g$ then 
$T_\lambda ^\sigma = \sigma T_\lambda \sigma ^{-1}$ is also a 
a conformal averaging operator, and 
\[
 T^\sigma _*(\mathfrak g) = \sigma (T_*(\mathfrak g)).
\]

Let $f,e,h$ be the standard basis of $\mathfrak{sl}_2$: 
$[h,e]=2e$, $[h,f]=-2f$, $[e,f]=h$.
For every nonzero $x\in \mathfrak{sl}_2$ there exists an automorphism
$\sigma $ such that either $\sigma (x)=\alpha h$ or $\sigma(x)=\alpha e$, 
$\alpha \in \mathbb C$, by the Jordan Normal Form Theorem applied 
to the 2-dimensional irreducible representation.

\begin{proposition}\label{prop:ConfAveSL_2}
 Let  $T_\lambda $ be a conformal averaging operator on $\mathfrak {sl}_2$. 
Then, up to an automorphism of $\mathfrak {sl}_2$, 
$T_\lambda $ is one of the following:
\begin{enumerate}
 \item $T_\lambda =\alpha \mathrm{id}$, $\alpha \in \mathbb C$;
 \item $T_\lambda(e)=T_\lambda (f)=0$, $T_\lambda (h)=\varphi(\lambda )h$; 
 \item $T_\lambda(e)=T_\lambda (h)=0$, $T_\lambda (f)=\varphi(\lambda )e$,
\end{enumerate}
where $\varphi(\lambda )\in \mathbb C[\lambda ]$ is an arbitrary nonzero polynomial.
\end{proposition}

\begin{proof}
Denote $L=T_*(\mathfrak {sl}_2)$. If $\dim L=3$
then $T_\lambda $ is non-degenerate.  By   Lemma~\ref{lem:Semisimple-NonDeg} $T_\lambda (x)=\alpha x$, $\alpha \in \mathbb C^*$.

Consider the case $\dim L=1$. Up to an automorphism, we may assume 
either $L=\mathbb Ch$ or $L=\mathbb Ce$.

In the first case, $T_\lambda $ is homogeneous, therefore, its 
structure is found in the previous section: $\Delta'=\varnothing $, $\mathfrak h_0^\perp =\mathbb Ch$, 
and thus we have $T_\lambda(e)=T_\lambda (f)=0$, $T_\lambda (h)=\varphi(\lambda )h$.

In the second case, $[T_\lambda (x),T_\mu(y)]=0$ for all $x,y\in \mathfrak {sl}_2$. Therefore, 
$T_\lambda ([e,\mathfrak {sl}_2])=\{0\}$, i.e., $e,h\in \Ker T_\lambda $. Hence, $T_\lambda (f)=\varphi(\lambda)e$ 
for some $\varphi (\lambda)\in \mathbb C[\lambda ]$.

Finally, assume $\dim L=2$. Such an operator may not be homogeneous, but 
it is easy to see that there are no 2-dimensional subalgebras in $\mathfrak{sl}_2$
that do not contain a Cartan subalgebra. Indeed, suppose $e\in L$ (up to automorphism),
and the second linear generator of $L$ is $x=\alpha f+\beta h$. Then $[e,x]\in L$ implies 
$h\in L$, which is impossible.
\end{proof}

Note that the same approach (classification of operators up to conjugation with an automorphism of $\mathfrak g$) works well 
for Rota---Baxter operators, i.e., linear maps $R:\mathfrak g\to \mathfrak g$ satisfying~\eqref{eq:RotaBaxer}.
In particular, for $\mathfrak g = \mathfrak{sl}_2$ the analogue of the proof of Proposition~\ref{prop:ConfAveSL_2} provides 
an easy description of Rota---Baxter operators of weight zero.

\begin{proposition}[c.f. {\cite{PeiBaiGuo2014}}]
Up to conjugation with an automorphism of $\mathfrak{g}=\mathfrak {sl}_2$ and scalar multiple, a Rota---Baxter operator of weight zero on
$\mathfrak {g}$ is one of the following:
\begin{enumerate}
 \item[\rm(R1)] $R=0$;
 \item[\rm(R2)] $R(e)=0$, $R(f)=t e - h$, $R(h)=2e$ ($t \in \mathbb C^*$);
 \item[\rm(R3)] $R(e)=0$, $R(f)=2te + h$, $R(h)=2e + \frac{1}{t}h$ ($t \in \mathbb C^*$);
 \item[\rm(R4)] $R(h)=h$, $R(e)=R(f)=0$;
 \item[\rm(R5)] $R(f)=h$, $R(e)=R(h)=0$;
 \item[\rm(R6)] $R(f)=e$, $R(e)=R(h)=0$.
\end{enumerate}
\end{proposition}

This description seems more compact than the one found in 
\cite{PeiBaiGuo2014}.

\begin{proof}
Indeed, suppose $R$ is a Rota---Baxter operator of  weight zero on $\mathfrak g = \mathfrak{sl}_2$, then $R(\mathfrak g)$ is a subalgebra 
of $\mathfrak g$ and $\Ker R$ is a module over $R(\mathfrak g)$. 
Nontrivial cases are: $\dim \Ker R=1,2$. 

If $\dim\Ker R=1$ then $\Ker R\subset R(\mathfrak g) $:  otherwise, $\mathfrak g = R(\mathfrak g)+\Ker R $  which implies 
$\Ker R$ to be an ideal of $\mathfrak g$. 
Up to conjugation with an automorphism, $\Ker R = \mathbb Ch$ or $\Ker R = \mathbb Ce$ (Jordan normal form of a matrix).

If $\Ker R =\mathbb Ch \subset R(\mathfrak g)$  then $R(\mathfrak g)$ is a 2-dimensional normalizer of the Cartan subalgebra, which is impossible. 

If $\Ker R = \mathbb Ce$ then $R(\mathfrak g) = \mathbb Ce +\mathbb Ch$ since it acts on $\Ker R$.
In this case, the operator $R$ is given in general by
\[
 R(e)=0,\quad R(f)=\alpha e+\beta h, \quad R(h)=\gamma e +\delta h,
\]
where the Rota---Baxter relation \eqref{eq:RotaBaxer} holds if and only if
\[
 \gamma \delta - 2\delta\beta =0, \quad \gamma^2-4\alpha \delta =-2\gamma \beta.
\]
If $\delta =0$ then, up to a scalar multiple, we have operator (R2). Analogously, if 
 $\delta \ne 0$ then we obtain (R3) [$\Ker R$ is 2-dimensional in the latter case].

Finally, suppose $\dim\Ker R=2$. Then $\dim R(\mathfrak g)=1 $ and, up to conjugation with an automorphism,
$ R(\mathfrak g)=\mathbb Ch$ or $ R(\mathfrak g)= \mathbb Ce$. 

If $ R(\mathfrak g)=\mathbb Ch$ then $\Ker R$  is a homogeneous subspace of $\mathfrak{g}$, so either 
$\Ker R=\mathbb C h+\mathbb C e$ (up to an involution) or 
$\Ker R=\mathbb C f+\mathbb C e$. These cases provide (R5) or (R4), respectively.

If $ R(\mathfrak g)=\mathbb Ce$ then $\Ker R=\mathbb C h+\mathbb C e$ is the only option 
for a 2-dimensional $e$-invariant subspace. Thus in this case we have (R6).

\end{proof}

\bibliography{jmp_average}

\providecommand{\noopsort}[1]{}\providecommand{\singleletter}[1]{#1}%
\begin{thebibliography}{10}

\bibitem{BaiGuoNi2012}
C.~Bai, O.~Bellier, L.~Guo, and X.~Ni.
\newblock Splitting of operations, {M}anin products, and {R}ota---{B}axter
  operators.
\newblock {\em Int. Math. Res. Notes}, 3:485--524, 2013.

\bibitem{BaiNiGuo2014}
C.~Bai, X.~Ni, and L.~Guo.
\newblock Generalizations of the classical {Y}ang---{B}axter equation and
  {O}-operators.
\newblock {\em J. Math. Phys.}, 52:063515, 2011.
\newblock DOI:10.1063/1.3600538.

\bibitem{BDK2001}
B.~Bakalov, A.~D'Andrea, and V.~G. Kac.
\newblock Theory of finite pseudoalgebras.
\newblock {\em Adv. Math.}, 162(1):1--140, 2001.

\bibitem{BelDr}
A.~A. Belavin and V.~G. Drinfeld.
\newblock The classical {Y}ang---{B}axter equation for simple {L}ie algebras.
\newblock {\em Funktsional. Anal. i Prilozhen.}, 17(3):69--70, 1983.

\bibitem{DK1998}
A.~D'Andrea and V.~G. Kac.
\newblock Structure theory of finite conformal algebras.
\newblock {\em Sel. Math., New Ser.}, 4:377--418, 1998.

\bibitem{FBZhvi2001}
E.~Frenkel and D.~Ben-Zvi.
\newblock {\em Vertex algebras and algebraic curves}, volume~88 of {\em
  Mathematical Surveys and Monographs}.
\newblock AMS, Providence, RI, second edition, 2004.

\bibitem{GubKol2011}
V.~Gubarev and P.~Kolesnikov.
\newblock Embedding of dendriform algebras into {R}ota---{B}axter algebras.
\newblock {\em Cent. Eur. J. Math.}, 11(2):226--245, 2013.

\bibitem{GubKol2014}
V.~Gubarev and P.~Kolesnikov.
\newblock Operads of decorated trees and their duals.
\newblock {\em Comment. Math. Univ. Carolin.}, 55(4):421--445, 2014.

\bibitem{GuoPei2014}
L.~Guo and J.~Pei.
\newblock Averaging algebras, {S}chr{\"o}der numbers and rooted trees.
\newblock {\em J. Algebraic Combinatorics}, 2014.
\newblock DOI:10.1007/s10801-014-0574-x.

\bibitem{Kac1996}
V.~G. Kac.
\newblock {\em Vertex algebras for beginners}, volume~10 of {\em University
  Lecture Series}.
\newblock AMS, Providence, RI, second edition, 1996.

\bibitem{Kol2015AL}
P.~S. Kolesnikov.
\newblock Homogeneous averaging operators on semisimple {L}ie algebras.
\newblock {\em Algebra and Logic}, 53(6):510--511, 2015.

\bibitem{MontStolinZelmanov}
F.~Montaner, A.~Stolin, and E.~Zelmanov.
\newblock Classification of {L}ie bialgebras over current algebras.
\newblock {\em Sel. Math. New Ser.}, 16:935--962, 2010.

\bibitem{PeiBaiGuo2014}
J.~Pei, C.~Bai, and L.~Guo.
\newblock Rota-baxter operators on {$sl(2,C)$} and solutions of the classical
  {Y}ang---{B}axter equation.
\newblock {\em J. Math. Phys.}, 55:021701, 2014.

\end{thebibliography}

\end{document}